\input amstex
\input amsppt.sty
\magnification=\magstep1
\hsize=33truecc
\vsize=22.2truecm
\baselineskip=16truept
\NoBlackBoxes
\TagsOnRight \pageno=1 \nologo
\def\Z{\Bbb Z}

\def\l{\left}
\def\r{\right}
\def\bg{\bigg}
\def\({\bg(}
\def\[{\bg\lfloor}
\def\){\bg)}
\def\]{\bg\rfloor}
\def\t{\text}
\def\f{\frac}

\def\p{\ (\roman{mod}\ p)}

\def\bi{\binom}
\def\eq{\equiv}

\def\ls{\leqslant}
\def\gs{\geqslant}
\def\mo{\roman{mod}}

\def\Proof{\noindent{\it Proof}}

\def\Remark{\medskip\noindent{\it  Remark}}

\def\Ack{\medskip\noindent {\bf Acknowledgment}}
\hbox {J. Number Theory 147(2015), no.\,1, 326--341.}
\bigskip
\topmatter
\title Supercongruences motivated by $e$\endtitle
\author Zhi-Wei Sun\endauthor
\leftheadtext{Zhi-Wei Sun}
\rightheadtext{Supercongruences motivated by $e$}
\affil Department of Mathematics, Nanjing University\\
 Nanjing 210093, People's Republic of China
  \\  zwsun\@nju.edu.cn
  \\ {\tt http://math.nju.edu.cn/$\sim$zwsun}
\endaffil
\abstract In this paper we establish some new supercongruences motivated by the well-known fact $\lim_{n\to\infty}(1+1/n)^n=e$. Let $p>3$ be a prime. We prove that
$$\sum_{k=0}^{p-1}\bi{-1/(p+1)}k^{p+1}\eq 0\ (\mo\ p^5)\ \, \t{and}\ \sum_{k=0}^{p-1}\bi{1/(p-1)}k^{p-1}\eq\f23p^4B_{p-3}\ (\mo\ p^5),$$
where $B_0,B_1,B_2,\ldots$ are Bernoulli numbers.
We also show that for any $a\in\Z$ with $p\nmid a$ we have
$$\sum_{k=1}^{p-1}\f1k\l(1+\f ak\r)^k\eq-1\pmod{p}\ \ \t{and}\ \ \sum_{k=1}^{p-1}\f1{k^2}\l(1+\f ak\r)^k\eq1+\f1{2a}\pmod{p}.$$
\endabstract
\thanks 2010 {\it Mathematics Subject Classification}.\,Primary 11B65;
Secondary 05A10,  11A07, 11B68.
\newline\indent Supported by the National Natural Science
Foundation (grant 11171140) of China and the PAPD of Jiangsu Higher Education Institutions.
\endthanks
\endtopmatter
\document

\heading{1. Introduction}\endheading

A $p$-adic congruence (with $p$ prime) is called a {\it
supercongruence} if it happens to hold modulo higher powers of $p$.
Here is a classical example due to J. Wolstenholme (cf. [W] or
[HT]):
$$\sum_{k=1}^{p-1}\f1k\eq0\ (\mo\ p^2)\ \ \t{and}\ \ \bi{2p-1}{p-1}\eq1\ (\mo\ p^3)$$
for every prime $p>3$. In 1900 Glaisher [G1, G2] showed further that
$$\bi{2p-1}{p-1}\eq1-\f23p^3B_{p-3}\pmod {p^4}$$
for any prime $p>3$, where $B_0,B_1,B_2,\ldots$ are Bernoulli numbers.
(See [IR, pp. 228--241] for an introduction to Bernoulli numbers.)
The reader may consult [L, Su1, Su3, T] for some other known supercongruences.

In this paper we establish some new supercongruences modulo prime
powers motivated by the well-known formula
$$\lim_{n\to\infty}\l(1+\f1n\r)^n=e.$$

Now we state our main results.

\proclaim{Theorem 1.1} Let $p>3$ be a prime. Then
$$\sum_{k=0}^{p-1}\bi{-1/(p+1)}k^{p+1}\eq0\ \ (\mo\ p^5).\tag1.1$$
Moreover, if $p>5$ then
$$\sum_{k=0}^{p-1}\bi{-1/(p+1)}k^{p+1}\eq\f{p^5}{18}B_{p-3}\ (\mo\ p^6).\tag1.2$$
\endproclaim

\proclaim{Theorem 1.2} Let $p>3$ be a prime and let $m$ be an integer not divisible by $p$. Then we have
$$\sum_{k=0}^{p-1}(-1)^{km}\bi{p/m-1}k^m\eq\f{(m-1)(7m-5)}{36m^2}p^4B_{p-3}\ \ (\mo\ p^5).\tag1.3$$
In particular,
$$\sum_{k=0}^{p-1}\bi{1/(p-1)}k^{p-1}\eq\f23 p^4B_{p-3}\ \ (\mo\ p^5).\tag1.4$$
\endproclaim

\Remark\ 1.1. Note that if $p$ is a prime and $m$ is an integer with $p\nmid m$ then
$$\bi{p/m-1}k\eq\bi{-1}k=(-1)^k\not\eq0\pmod p\quad\t{for all}\ k=0,1,\ldots,p-1.$$
(1.1)-(1.4) are interesting since supercongruences modulo $p^5$ are very rare.
We conjecture that there are no composite numbers $p>1$ satisfying (1.1) or the congruence
$$\sum_{k=0}^{p-1}\bi{1/(p-1)}k^{p-1}\eq0\pmod {p^4}.\tag1.5$$

\medskip
\proclaim{Theorem 1.3} Let $p>3$ be a prime and let $m$ be an integer not divisible by $p$.

{\rm (i)} If $p>5$, then
$$\sum_{k=1}^{p-1}\f{(-1)^{km}}{k^2}\bi{p/m-1}k^m\eq\f 1p\sum_{k=1}^{p-1}\f1k\ \ (\mo\ p^3).\tag1.6$$
Also, for any $n=1,\ldots,(p-3)/2$ we have
$$\sum_{k=1}^{p-1}\f{(-1)^{km}}{k^{2n}}\bi{p/m-1}k^m\eq -\f p{2n+1}B_{p-1-2n}\ \ (\mo\ p^2).\tag1.7$$

{\rm (ii)} For every $n=1,\ldots,(p-3)/2$, we have
$$\sum_{k=1}^{p-1}\f{(-1)^{km}}{k^{2n-1}}\bi{p/m-1}k^m\eq \l(1+\f{1-m}{2m}(2n+1)\r)\f {p^2n}{2n+1}B_{p-1-2n}\ (\mo\ p^3).\tag1.8$$
\endproclaim
\Remark\ 1.2. If $n$ is a positive integer and $p>2n+1$ is a prime, then (1.8) with $m=p\pm1$ yields the congruences
$$\sum_{k=1}^{p-1}\f1{k^{2n-1}}\bi{1/(p-1)}k^{p-1}\eq-\f{2p^2n^2}{2n+1}B_{p-1-2n}\ \ (\mo\ p^3)\tag1.9$$
and
$$\sum_{k=1}^{p-1}\f1{k^{2n-1}}\bi{-1/(p+1)}k^{p+1}\eq\f{p^2n}{2n+1}B_{p-1-2n}\ \ (\mo\ p^3).\tag1.10$$

\proclaim{Theorem 1.4} Let $p$ be an odd prime and let $a\in\Z$ with $p\nmid a$. Then
$$\sum_{k=1}^{p-1}\f1k\l(1+\f ak\r)^k\eq-1\pmod p.\tag1.11$$
If $p>3$, then
$$\sum_{k=1}^{p-1}\f1{k^2}\l(1+\f ak\r)^k\eq1+\f1{2a}\pmod p.\tag1.12$$
\endproclaim
\Remark\ 1.3. (1.11) with $a=1$ yields the congruence
$$\sum_{k=1}^{p-1}\f{(k+1)^k}{k^{k+1}}\eq-1\pmod p.\tag1.13$$

\medskip
We will show Theorems 1.1-1.2 in the next section. Theorems 1.3 and 1.4 will be proved in Sections 3 and 4 respectively.

To conclude this section, we pose two related conjectures for further research.

 \proclaim{Conjecture 1.1} Let $m>2$ and $q>0$ be integers with $m$ even or $q$ odd. Then, for any prime $p>mq$ we have the supercongruence
 $$\sum_{k=0}^{p-1}(-1)^{km}\bi{p/m-q}k^m\eq0\ \ (\mo\ p^3).\tag1.14$$
 \endproclaim
 \Remark\ 1.4. Clearly (1.14) with $q=1$ follows from (1.3).
 \medskip

For a prime $p$ and a $p$-adic number $x$, as usual we let $\nu_p(x)$ denote
the $p$-adic valuation (i.e., $p$-adic order) of $x$.

 \proclaim{Conjecture 1.2} Let $p$ be a prime and let $n$ be a
 positive integer. Then
 $$\nu_p\(\sum_{k=0}^{n-1}\bi{-1/(p+1)}k^{p+1}\)\gs c_p\l\lfloor\f{\nu_p(n)+1}2\r\rfloor,\tag1.15$$
 where
 $$c_p=\cases1&\t{if}\ p=2,\\3&\t{if}\ p=3,\\5&\t{if}\ p\gs5.\endcases$$
 If $p>3$, then we have
$$\nu_p\(\sum_{k=0}^{n-1}\bi{1/(p-1)}k^{p-1}\)\gs 4\l\lfloor\f{\nu_p(n)+1}2\r\rfloor.\tag1.16$$
 \endproclaim

\heading{2. Proofs of Theorems 1.1 and 1.2}\endheading

For $m=1,2,3,\ldots$ and $n=0,1,2,\ldots$, we
define
$$H_n^{(m)}:=\sum_{0<k\ls n}\f1{k^m}$$
and call it a harmonic number of order $m$. Those $H_n=H_n^{(1)}\
(n=0,1,2,\ldots)$ are usually called {\it harmonic numbers}.

\proclaim{Lemma 2.1} Let $p>3$ be a prime. Then
$$H_{p-1}\eq-\f{p^2}3B_{p-3}\ (\mo\ p^3)\ \ \t{and}\ \ \ H_{p-1}^{(2)}\eq \f23pB_{p-3}\ (\mo\ p^2).\tag2.1$$
Also,
$$\sum_{k=1}^{p-1}H_k^{(2)}\eq p^2B_{p-3}\ (\mo\ p^3),\ \
\sum_{k=1}^{p-1}H_k^{(3)}\eq-\f23pB_{p-3}\ (\mo\ p^2),\tag2.2$$
and
$$\sum_{k=1}^{p-1}H_k^{(4)}\eq H_{p-1}^{(3)}\eq0\pmod p.\tag2.3$$
\endproclaim
\Proof. The two congruences in (2.1) are known results due to Glaisher [G2], see also Theorem 5.1 and Corollary 5.1 of [S].

For $m=2,3,4$ we have
$$\sum_{k=1}^{p-1}H_k^{(m)}=\sum_{k=1}^{p-1}\sum_{j=1}^k\f1{j^m}
=\sum_{j=1}^{p-1}\f{\sum_{k=j}^{p-1}1}{j^m}=\sum_{j=1}^{p-1}\f{p-j}{j^m}=pH_{p-1}^{(m)}-H_{p-1}^{(m-1)}.$$
Note also that
$$H_{p-1}^{(3)}=\sum_{k=1}^{(p-1)/2}\(\f1{k^3}+\f1{(p-k)^3}\)\eq0\pmod p.$$
Combining these with (2.1), we immediately get (2.2) and (2.3). \qed

\proclaim{Lemma 2.2} Let $p>3$ be a prime. Set
$$\Sigma_1=\sum_{k=1}^{p-1}\sum_{1\ls i<j\ls k}\l(\f1{ij^2}+\f1{i^2j}\r),\ \
\Sigma_2=\sum_{k=1}^{p-1}\sum_{1\ls i<j\ls k}\f1{i^2j^2},$$
$$\Sigma_3=\sum_{k=1}^{p-1}\sum_{1\ls i<j\ls k}\l(\f1{ij^3}+\f1{i^3j}\r)
\ \ \t{and}\ \ \ \Sigma_4=\sum_{k=1}^{p-1}\sum_{1\ls i<j\ls k}\f{H_k^{(2)}}{ij}.$$
Then we have
$$\Sigma_1\eq pB_{p-3}\pmod{p^2},\qquad \ \Sigma_2\eq B_{p-3}\pmod p\tag2.4$$
and
$$\Sigma_3\eq\Sigma_4\eq-B_{p-3}\pmod p.\tag2.5$$
\endproclaim
\Proof. With the help of Lemma 2.1,
$$\align\Sigma_1=&\sum_{1\ls i<j\ls p-1}\l(\f1{ij^2}+\f1{i^2j}\r)\sum_{k=j}^{p-1}1
\\=&p\sum_{1\ls i<j\ls p-1}\l(\f1{ij^2}+\f1{i^2j}\r)-\sum_{1\ls i<j\ls p-1}\l(\f1{ij}+\f1{i^2}\r)
\\=&p\l(H_{p-1}H_{p-1}^{(2)}-H_{p-1}^{(3)}\r)-\sum_{1\ls i<j\ls p-1}\f1{ij}-\sum_{i=1}^{p-1}\f{p-1-i}{i^2}
\\\eq&0-\f12\l(H_{p-1}^2-H_{p-1}^{(2)}\r)-(p-1)H_{p-1}^{(2)}+H_{p-1}\eq pB_{p-3}\pmod{p^2}.
\endalign$$
Recall the congruence $\sum_{k=1}^{p-1}H_k/k^2\eq B_{p-3}\pmod p$ (cf. [ST, (5.4)]). Note that
$$\align \Sigma_2=&\sum_{1\ls i<j\ls p-1}\f1{i^2j^2}\sum_{k=j}^{p-1}1=p\sum_{1\ls i<j\ls p-1}\f1{i^2j^2}-\sum_{1\ls i<j\ls p-1}\f1{i^2j}
\\\eq&-\sum_{i=1}^{p-1}\f{H_{p-1}-H_i}{i^2}\eq\sum_{k=1}^{p-1}\f{H_k}{k^2}\eq B_{p-3}\pmod p
\endalign$$
and
$$\align \Sigma_3=&\sum_{1\ls i<j\ls p-1}\l(\f1{ij^3}+\f1{i^3j}\r)\sum_{k=j}^{p-1}1=\sum_{1\ls i<j\ls p-1}\l(\f1{ij^3}+\f1{i^3j}\r)(p-j)
\\\eq&-\sum_{1\ls i<j\ls p-1}\l(\f1{ij^2}+\f1{i^3}\r)=-\sum_{j=1}^{p-1}\f{H_j-1/j}{j^2}-\sum_{i=1}^{p-1}\f{p-1-i}{i^3}
\\\eq&-\sum_{k=1}^{p-1}\f{H_k}{k^2}+2H_{p-1}^{(3)}+H_{p-1}^{(2)}\eq-B_{p-3}\pmod p.
\endalign$$

 In light of (2.2),
$$\align \Sigma_4=&\sum_{1\ls i<j\ls p-1}\f1{ij}\(\sum_{k=1}^{p-1}H_k^{(2)}-\sum_{s=1}^{j-1}H_s^{(2)}\)
\\\eq&-\sum_{1\ls i<j\ls p-1}\f1{ij}\sum_{1\ls t\ls s<j}\f1{t^2}=-\sum_{1\ls i<j\ls p-1}\f1{ij}\sum_{t=1}^j\f{j-t}{t^2}
\\=&-\sum_{j=1}^{p-1}H_{j-1}H_j^{(2)}+\sum_{j=1}^{p-1}\f{H_j}j\l(H_j-\f1j\r)\pmod {p^2}.
\endalign$$
For every $k=1,\ldots,p-1$, we have
$$H_{p-k}=H_{p-1}-\sum_{0<j<k}\f1{p-j}\eq H_{k-1}\pmod p$$
and
$$H_{p-k}^{(2)}=H_{p-1}^{(2)}-\sum_{0<j<k}\f1{(p-j)^2}\eq-H_{k-1}^{(2)}\pmod p.$$
Thus
$$\align \sum_{k=1}^{p-1}H_{k-1}H_k^{(2)}\eq& \sum_{k=1}^{p-1}H_{p-k}H_k^{(2)}=\sum_{k=1}^{p-1}H_kH_{p-k}^{(2)}
\\\eq&-\sum_{k=1}^{p-1}H_kH_{k-1}^{(2)}=-\sum_{k=1}^{p-1}H_k\l(H_k^{(2)}-\f1{k^2}\r)\pmod p.
\endalign$$
and hence
$$\Sigma_4\eq\sum_{k=1}^{p-1}H_kH_k^{(2)}-2\sum_{k=1}^{p-1}\f{H_k}{k^2}+\sum_{k=1}^{p-1}\f{H_k^2}k\pmod p.\tag2.6$$
Since
$$\sum_{k=1}^{p-1}\f{H_k^2}k=\sum_{k=1}^{p-1}\f{H_{p-k}^2}{p-k}\eq-\sum_{k=1}^{p-1}\f1k\l(H_k-\f1k\r)^2=-\sum_{k=1}^{p-1}\f{H_k^2}k+2\sum_{k=1}^{p-1}\f{H_k}{k^2}-H_{p-1}^{(3)}
\ (\mo\ p)$$
and $H_{p-1}^{(3)}\eq0\pmod p$, we have
$$\sum_{k=1}^{p-1}\f{H_k^2}k\eq\sum_{k=1}^{p-1}\f{H_k}{k^2}\eq B_{p-3}\pmod p.\tag2.7$$
Combining this with (2.6) and the congruence $\sum_{k=1}^{p-1}H_kH_k^{(2)}\eq0\pmod p$ (cf. [Su2, (2.8)]), we obtain
$\Sigma_4\eq -B_{p-3}\pmod p$. This concludes the proof. \qed

\medskip
\noindent{\it Proof of Theorem} 1.1. (1.1) in the case $p=5$ can be verified directly; in fact,
$$\sum_{k=0}^{5-1}\bi{-1/(5+1)}k^{5+1}\eq 5^5\pmod{5^6}.$$
Below we assume $p>5$.
As (1.2) implies (1.1), we only need to prove (1.2).

For any $p$-adic integer $x$, we can write $(1+px)^m$ as the $p$-adic series $\sum_{n=0}^\infty\bi mnp^nx^n$.
Thus, for each $k=1,\ldots,p-1$, we have
$$\align &\bi{-1/(p+1)}k^{p+1}=\bi{p/(p+1)-1}k^{p+1}=\prod_{j=1}^k\l(1-\f
p{(p+1)j}\r)^{p+1}
\\\eq&\prod_{j=1}^k\l(1-\f{p}{j}+\bi{p+1}2\f{p^2}{(p+1)^2j^2}
-\bi{p+1}3\f{p^3}{(p+1)^3j^3}+\bi{p+1}4\f{p^4}{(p+1)^4j^4}\r)
\\=&\prod_{j=1}^k\l(1-\f pj+\f{p^3}{2(p+1)j^2}-\f{p^4(p-1)}{6(p+1)^2j^3}+\f{p^5(p-1)(p-2)}{24(p+1)^3j^4}\r)
\\\eq&\prod_{j=1}^k\l(1-\f pj+\f{p^3(p^2-p+1)}{2j^2}-\f{p^4}{6j^3}(p-1)(1-2p)+\f{2p^5}{24j^4}\r)
\\\eq&\prod_{j=1}^k\l(1-\f pj+\f{p^5-p^4+p^3}{2j^2}+\f{p^4-3p^5}{6j^3}+\f{p^5}{12j^4}\r)\pmod{p^6}
\endalign$$
and hence
$$\align&\bi{-1/(p+1)}k^{p+1}-\prod_{j=1}^k\l(1-\f pj\r)
\\\eq&\f{p^5-p^4+p^3}2H_k^{(2)}+\f{p^4-3p^5}6H_k^{(3)}+\f{p^5}{12}H_k^{(4)}
\\&+\f {p(p^4-p^3)}2\sum_{1\ls i<j\ls k}\l(\f1{ij^2}+\f1{i^2j}\r)-\f{p^5}6\sum_{1\ls i<j\ls k}\l(\f1{ij^3}+\f1{i^3j}\r)
\\&+\f{p^5}2\sum_{1\ls i_1<i_2\ls k}\f1{i_1i_2}\(\sum_{j=1}^k\f1{j^2}-\f1{i_1^2}-\f1{i_2^2}\)\ (\mo\ p^6).
\endalign$$
Thus, in view of Lemmas 2.1 and 2.2, we obtain
$$\align&\sum_{k=1}^{p-1}\bi{-1/(p+1)}k^{p+1}-\sum_{k=1}^{p-1}(-1)^k\bi{p-1}k
\\\eq&\f{p^5-p^4+p^3}2p^2B_{p-3}+\f{p^4-3p^5}6\l(-\f23pB_{p-3}\r)+\f{p^5}{12}\times0
\\&+\f{p^5-p^4}2pB_{p-3}-\f{p^5}6(-B_{p-3})+\f{p^5}2\l(-B_{p-3}-(-B_{p-3})\r)
\\\eq&\f{p^5}{18}B_{p-3}\pmod{p^6}
\endalign$$
and hence (1.2) follows since $\sum_{k=0}^{p-1}(-1)^k\bi{p-1}k=(1-1)^{p-1}=0$.

The proof of Theorem 1.1 is now complete. \qed

\medskip
\noindent{\it Proof of Theorem} 1.2. For each $k\in\{1,\ldots,p-1\}$, obviously
$$(-1)^{km}\bi{p/m-1}k^m=\prod_{j=1}^k\l(1-\f{p}{jm}\r)^m$$
is congruent to
$$\align&\prod_{j=1}^k\l(1-\f{pm}{jm}+\f{m(m-1)}2\cdot\f{p^2}{j^2m^2}-\bi m3\f{p^3}{j^3m^3}+\bi m4\f{p^4}{j^4m^4}\r)
\\=&\prod_{j=1}^k\l(1-\f pj+\f{m-1}{2m}\cdot\f{p^2}{j^2}-\f{(m-1)(m-2)}{6m^2}\cdot\f{p^3}{j^3}+\f{(m-1)(m-2)(m-3)}{24m^3}\cdot\f{p^4}{j^4}\r)
\endalign$$ modulo $p^5$, and hence
$$\align &(-1)^{km}\bi{p/m-1}k^m-\prod_{j=1}^k\l(1-\f pj\r)
\\\eq&\f{m-1}{2m}p^2H_k^{(2)}-\f{(m-1)(m-2)}{6m^2}p^3H_k^{(3)}+\f{(m-1)(m-2)(m-3)}{24m^3}p^4H_k^{(4)}
\\&-\f{m-1}{2m}p^3\sum_{1\ls i<j\ls k}\l(\f1{ij^2}+\f1{i^2j}\r)+\f{(m-1)^2}{4m^2}p^4\sum_{1\ls i<j\ls k}\f1{i^2j^2}
\\&+\f{(m-1)(m-2)}{6m^2}p^4\sum_{1\ls i<j\ls k}\l(\f1{ij^3}+\f1{i^3j}\r)
\\&+\f{m-1}{2m}p^4\sum_{1\ls i_1<i_2\ls k}\f1{i_1i_2}\sum^k\Sb j=1\\j\not=i_1,i_2\endSb\f1{j^2}
\ \ (\mo\ p^5).
\endalign$$
Therefore, applying (2.2) and (2.3) we get
$$\align&\sum_{k=1}^{p-1}(-1)^{km}\bi{p/m-1}k^m-\sum_{k=1}^{p-1}(-1)^k\bi{p-1}k
\\\eq&\f{m-1}{2m}p^2(p^2B_{p-3})-\f{(m-1)(m-2)}{6m^2}p^3\l(-\f23pB_{p-3}\r)
\\&-\f{m-1}{2m}p^3\Sigma_1+\f{(m-1)^2}{4m^2}p^4\Sigma_2
\\&+\f{(m-1)(m-2)}{6m^2}p^4\Sigma_3+\f{m-1}{2m}p^4(\Sigma_4-\Sigma_3)\pmod{p^5},
\endalign$$
where $\Sigma_1,\Sigma_2,\Sigma_3,\Sigma_4$ are defined in Lemma 2.2.
Combining this with Lemma 2.2, we finally obtain
$$\align &\sum_{k=0}^{p-1}(-1)^{km}\bi{p/m-1}k^m-(1-1)^{p-1}
\\\eq&\f{m-1}{2m}p^4B_{p-3}+\f{(m-1)(m-2)}{9m^2}p^4B_{p-3}
\\&-\f{m-1}{2m}p^3(pB_{p-3})+\f{(m-1)^2}{4m^2}p^4B_{p-3}+\f{(m-1)(m-2)}{6m^2}p^4(-B_{p-3})
\\=&\l(\f{(m-1)^2}{4m^2}-\f{(m-1)(m-2)}{18m^2}\r)p^4B_{p-3}\pmod {p^5},
\endalign$$
which gives (1.3). Putting $m=p-1$ in (1.3) we immediately get (1.4). This concludes the proof. \qed

\heading{3. Proof of Theorem 1.3}\endheading

\proclaim{Lemma 3.1} Let $m$ and $n$ be positive integers with $m\ls 2n$,
and let $p>2n+1$ be a prime. Then
$$\sum_{k=1}^{p-1}\f{H_{k}^{(m)}}{k^{2n+1-m}}\eq\f{(-1)^{m-1}}{2n+1}\bi{2n+1}m B_{p-1-2n}\ (\mo\ p).\tag 3.1$$
When $m<2n$ we have
$$\sum_{k=1}^{p-1}\f{H_{k}^{(m)}}{k^{2n-m}}\eq\f{pB_{p-1-2n}}{2n+1}\l(n+(-1)^m\f{n-m}{m+1}\bi{2n+1}m\r)\ (\mo\ p^2).\tag 3.2$$
\endproclaim
\Proof. Since $\sum_{k=1}^{p-1}k^s\eq0\ (\mo\ p)$ for any integer
$s\not\eq0\ (\mo\ p-1)$ (see, e.g., [IR, p.\,235]),
 we have $\sum_{k=1}^{p-1}1/k^{2n+1}\eq0\ (\mo\ p)$. Hence
$$\align &\sum_{k=1}^{p-1}\f{H_{k}^{(m)}}{k^{2n+1-m}}
\\\eq&\sum_{k=1}^{p-1}\f1{k^{2n+1}}+\sum_{k=1}^{p-1}\f1{k^{2n+1-m}}\sum_{0<j<k}j^{p-1-m}
\\\eq&\sum_{k=1}^{p-1}\f1{k^{2n+1-m}(p-m)}\sum_{i=0}^{p-1-m}\bi{p-m}iB_ik^{p-m-i}\ \ \ \ (\t{by [IR, p.\,230]})
\\\eq&-\f1m\sum_{i=0}^{p-1-m}\bi{p-m}iB_i\sum_{k=1}^{p-1}k^{p-1-2n-i}
\\\eq&\f1m\sum\Sb 0\ls i\ls p-1-m\\p-1\mid 2n+i\endSb\bi{p-m}iB_i=\f1m\bi{p-m}{2n+1-m}B_{p-1-2n}
\\\eq&\f1m\bi{-m}{2n+1-m}B_{p-1-2n}=\f{(-1)^{m-1}}{2n+1}\bi{2n+1}mB_{p-1-2n}\pmod p.
\endalign$$
This proves (3.1).

 Now assume that $m<2n$. As $m,2n-m\in\{1,\ldots,p-2\}$, we have
 $$\sum_{j=1}^{p-1}\f1{j^m}\eq\sum_{k=1}^{p-1}\f1{k^{2n-m}}\eq0\pmod{p}.$$
 It is known that
$$\sum_{k=1}^{p-1}\f1{k^s}\eq\f{ps}{s+1}B_{p-1-s}\ (\mo\ p^2)\quad\t{for each}\ s=1,\ldots,p-2.\tag3.3$$
(See, e.g., [G2] or [S, Corollary 5.1].) Thus
$$\sum_{j=1}^{p-1}\f1{j^m}\sum_{k=1}^{p-1}\f1{k^{2n-m}}+\sum_{k=1}^{p-1}\f1{k^{2n}}
\eq\sum_{k=1}^{p-1}\f1{k^{2n}}\eq\f{2n}{2n+1}pB_{p-1-2n}\pmod{p^2}.$$
On the other hand,
$$\align&\sum_{j=1}^{p-1}\f1{j^m}\sum_{k=1}^{p-1}\f1{k^{2n-m}}+\sum_{k=1}^{p-1}\f1{k^{2n}}-\sum_{1\ls j\ls k\ls p-1}\f1{j^mk^{2n-m}}
\\=&\sum_{1\ls k\ls j\ls p-1}\f1{j^mk^{2n-m}}=\sum_{1\ls j\ls k\ls p-1}\f1{(p-j)^m(p-k)^{2n-m}}
\\=&\sum_{1\ls j\ls k\ls p-1}\f{(p+j)^m(p+k)^{2n-m}}{(p^2-j^2)^m(p^2-k^2)^{2n-m}}
\\\eq&\sum_{1\ls j\ls k\ls p-1}\f{(j^m+pmj^{m-1})(k^{2n-m}+p(2n-m)k^{2n-m-1})}{j^{2m}k^{2(2n-m)}}
\\\eq&\sum_{1\ls j\ls k\ls p-1}\l(\f1{j^mk^{2n-m}}+\f{pm}{j^{m+1}k^{2n-m}}+\f{p(2n-m)}{j^mk^{2n-m+1}}\r)\pmod{p^2}.
\endalign$$
Therefore, with the help of (3.1) we have
$$\align &\f{2n}{2n+1}pB_{p-1-2n}-2\sum_{k=1}^{p-1}\f{H_{k}^{(m)}}{k^{2n-m}}
\\\eq&pm\sum_{k=1}^{p-1}\f{H_{k}^{(m+1)}}{k^{2n-m}}+p(2n-m)\sum_{k=1}^{p-1}\f{H_{k}^{(m)}}{k^{2n-m+1}}
\\\eq&pm\f{(-1)^m}{2n+1}\bi{2n+1}{m+1}B_{p-1-2n}+p(2n-m)\f{(-1)^{m-1}}{2n+1}\bi{2n+1}mB_{p-1-2n}
\\=&(-1)^m\f{2(m-n)}{m+1}\bi{2n+1}m\f{pB_{p-1-2n}}{2n+1}\pmod{p^2}
\endalign$$
and hence (3.2) holds. \qed
\Remark\ 3.1. By [ST, (5.4)], $\sum_{k=1}^{p-1}H_k/k^2\eq B_{p-3}\ (\mo\ p)$ for any prime $p>3$.
By [M, (5)], $\sum_{k=1}^{p-1}H_k/k^3\eq-pB_{p-5}/10\ (\mo\ p^2)$ for any prime $p>5$.
Obviously these two results are particular cases of Lemma 3.1.

\proclaim{Lemma 3.2} Let $p>5$ be a prime. Then
$$\sum_{k=1}^{p-1}\f{1-pH_k}{k^2}\eq\f{H_{p-1}}p\ (\mo\ p^3).\tag3.4$$
\endproclaim
\Proof. In view of Theorem 5.1(a) and Remark 5.1 of [S],
$$\f{H_{p-1}^{(2)}}2\eq p\l(\f{B_{2p-4}}{2p-4}-2\f{B_{p-3}}{p-3}\r)\eq-\f{H_{p-1}}p\ (\mo\ p^3)$$
and $H_{p-1}^{(3)}\eq0\ (\mo\ p^2)$. Also,
$$\sum_{k=1}^{p-1}\f{H_{k-1}}{k^2}=\sum_{1\ls j<k\ls p-1}\f1{jk^2}\eq-3\f{H_{p-1}}{p^2}\ (\mo\ p^2)$$
by [T, Theorem 2.3]. So we have
$$\sum_{k=1}^{p-1}\f{1-pH_k}{k^2}=H_{p-1}^{(2)}-pH_{p-1}^{(3)}-p\sum_{k=1}^{p-1}\f{H_{k-1}}{k^2}
\eq\f{H_{p-1}}p\ (\mo\ p^3).$$
This concludes the proof. \qed

\medskip
\noindent{\it Proof of Theorem} 1.3. Let $k\in\{1,\ldots,p-1\}$. Then
$$\align(-1)^{km}\bi{p/m-1}k^m=&\prod_{j=1}^k\l(1-\f p{jm}\r)^m
\\\eq&\prod_{j=1}^k\l(1-\f pj+\f{m(m-1)}2\cdot\f{p^2}{j^2m^2}\r)
\\\eq&1-pH_k+\f{m-1}{2m}p^2H_k^{(2)}+p^2\sum_{1\ls i<j\ls k}\f1{ij}\pmod{p^3}.
\endalign$$
Thus, for every $r=1,\ldots,p-3$ we have
$$\aligned &\sum_{k=1}^{p-1}\f{(-1)^{km}}{k^r}\bi{p/m-1}k^m
\\\eq&\sum_{k=1}^{p-1}\f{1-pH_k}{k^r}+\f{m-1}{2m}p^2\sum_{k=1}^{p-1}\f{H_k^{(2)}}{k^r}+\f{p^2}2\sum_{k=1}^{p-1}\f{H_k^2-H_k^{(2)}}{k^r}\ (\mo\ p^3).
\endaligned\tag3.5$$

(i) Now suppose that $n\in\{1,\ldots,(p-3)/2\}$. Then (3.5) with $r=2n$ yields
$$\sum_{k=1}^{p-1}\f{(-1)^{km}}{k^{2n}}\bi{p/m-1}k^m\eq\sum_{k=1}^{p-1}\f{1-pH_k}{k^{2n}}\pmod{p^2}.\tag3.6$$
By (3.3) and Lemma 3.1,
$$\sum_{k=1}^{p-1}\f{1-pH_k}{k^{2n}}\eq\l(\f{2n}{2n+1}-1\r)pB_{p-1-2n}=-\f{pB_{p-1-2n}}{2n+1}\pmod{p^2}.$$
So (1.7) follows from (3.6).
If $n<(p-3)/2$, then
$$2\sum_{k=1}^{p-1}\f{H_k^{(2)}}{k^{2n}}\eq\sum_{k=1}^{p-1}\l(\f{H_k^{(2)}}{k^{2n}}+\f{H_{p-k}^{(2)}}{(p-k)^{2n}}\r)
\eq\sum_{k=1}^{p-1}\f{1/k^2}{k^{2n}}\eq0\pmod p,$$
and hence by (3.5) with $r=2n$ we obtain
$$\sum_{k=1}^{p-1}\f{(-1)^{km}}{k^{2n}}\bi{p/m-1}k^m\eq\sum_{k=1}^{p-1}\f{1-pH_k}{k^{2n}}+\f{p^2}2\sum_{k=1}^{p-1}\f{H_k^2}{k^{2n}}
\ (\mo\ p^3).\tag3.7$$
When $p>5$, (3.7) in the case $n=1$, together with (3.4) and the subtle congruence
$$\sum_{k=1}^{p-1}\f{H_k^2}{k^2}\eq0\ (\mo\ p)$$
of Sun [Su2, (1.5)], yields (1.6).

(ii) Fix $n\in\{1,\ldots,(p-3)/2\}$. Putting $r=2n-1$ in (3.5) we get
$$\aligned &\sum_{k=1}^{p-1}\f{(-1)^{km}}{k^{2n-1}}\bi{p/m-1}k^m
\\\eq&\sum_{k=1}^{p-1}\f{1-pH_k}{k^{2n-1}}+\f{m-1}{2m}p^2\sum_{k=1}^{p-1}\f{H_k^{(2)}}{k^{2n-1}}
+\f{p^2}2\sum_{k=1}^{p-1}\f{H_k^2-H_k^{(2)}}{k^{2n-1}}\ (\mo\ p^3).
\endaligned\tag3.8$$
In view of a known result (cf. [G2] or [S, Theorem 5.1(a)]),
$$\sum_{k=1}^{p-1}\f1{k^{2n-1}}\eq\f{n-2n^2}{2n+1}p^2B_{p-1-2n}\pmod{p^3}.$$
By Lemma 3.1,
$$\sum_{k=1}^{p-1}\f{H_k}{k^{2n-1}}\eq\f{1+3n-2n^2}{2(2n+1)}pB_{p-1-2n}\pmod{p^2}$$
and $$\sum_{k=1}^{p-1}\f{H_k^{(2)}}{k^{2n-1}}\eq-nB_{p-1-2n}\pmod{p}.$$
Note also that
$$\sum_{k=1}^{p-1}\f{H_k^2}{k^{2n-1}}=\sum_{k=1}^{p-1}\f{H_{p-k}^2}{(p-k)^{2n-1}}
\eq-\sum_{k=1}^{p-1}\f{(H_k-1/k)^2}{k^{2n-1}}\pmod p$$
and hence
$$\sum_{k=1}^{p-1}\f{H_k^2}{k^{2n-1}}\eq\sum_{k=1}^{p-1}\f{H_k}{k^{2n}}-\f12\sum_{k=1}^{p-1}\f1{k^{2n+1}}
\eq\sum_{k=1}^{p-1}\f{H_k}{k^{2n}}\eq B_{p-1-2n}\pmod p$$
with the help of (3.1) in the case $m=1$.
Combining all these we obtain (1.8) from (3.8).

The proof of Theorem 1.3 is now complete. \qed

\heading{4. Proof of Theorem 1.4}\endheading

\proclaim{Lemma 4.1} Let $p$ be an odd prime. Then, for any positive integers $d$ and $r$ with $d+r<p$, we have
$$\sum_{k=r}^{p-1}\f{\bi kr}{k^{r+d}}\eq0\pmod p.\tag4.1$$
\endproclaim
\Proof. Observe that
$$\align \sum_{k=r}^{p-1}\bi krk^{-r-d}=&\sum_{s=0}^{p-1-r}\bi{r+s}s(r+s)^{-r-d}=\sum_{s=0}^{p-1-r}(-1)^s\bi{-r-1}s(r+s)^{-r-d}
\\\eq&\sum_{s=0}^{p-1-r}\bi{p-1-r}s(-1)^s(s-1-(p-1-r))^{p-1-r-d}
\\=&\sum_{k=0}^{p-1-r}\bi{p-1-r}k(-1)^{p-1-r-k}(-1-k)^{p-1-r-d}
\\=&(-1)^d\sum_{k=0}^{p-1-r}\bi{p-1-r}k(-1)^k(k+1)^{p-1-r-d}\pmod p.
\endalign$$
It is known that for any positive integer $n$ we have
$$\sum_{k=0}^n\bi nk(-1)^k k^m=0\quad\t{for all}\ m=0,1,\ldots,n-1.$$
(See, e,g., [vLW, pp.\,125-126].) Therefore, (4.1) follows from the above. \qed

\medskip
\noindent{\it Proof of Theorem} 1.4. Let $d\in\{1,2\}$. With the help of Lemma 4.1, we have
$$\align\sum_{k=1}^{p-1}\f1{k^d}\l(1+\f ak\r)^k=&\sum_{k=1}^{p-1}\f{(k+a)^k}{k^{k+d}}
\\=&\sum_{k=1}^{p-1}\f1{k^{k+d}}\(k^k+\sum_{r=1}^k\bi krk^{k-r}a^r\)
\\=&\sum_{k=1}^{p-1}\f1{k^d}+\sum_{r=1}^{p-1}a^r\sum_{k=r}^{p-1}\f{\bi kr}{k^{r+d}}
\\\eq&\sum_{k=1}^{p-1}\f1{k^d}+\sum_{r=p-d}^{p-1}a^r\sum_{k=r}^{p-1}\f{\bi kr}{k^{r+d}}\ \pmod p.
\endalign$$
Thus,
$$\sum_{k=1}^{p-1}\f1k\l(1+\f ak\r)^k\eq\sum_{k=1}^{(p-1)/2}\l(\f1k+\f1{p-k}\r)+\f{a^{p-1}}{(p-1)^p}\eq-1\pmod p$$
in view of Fermat's little theorem. When $p>3$, we have $\sum_{k=1}^{p-1}1/k^2\eq0\pmod p$ by [W], and hence
$$\align\sum_{k=1}^{p-1}\f1{k^2}\l(1+\f ak\r)^k\eq&\sum_{k=1}^{p-1}\f1{k^2}+\f{a^{p-1}}{(p-1)^{p+1}}+a^{p-2}\(\f{\bi{p-2}{p-2}}{(p-2)^p}+\f{\bi{p-1}{p-2}}{(p-1)^p}\)
\\\eq&0+1+\f1a\l(\f1{-2}+1\r)=1+\f1{2a}\ \pmod p
\endalign$$
as desired. This concludes the proof. \qed
\medskip

\Ack. The author would like to thank the referee for helpful comments.


\medskip

 \widestnumber\key{vLW}

 \Refs

\ref\key G1\by J.W.L. Glaisher\paper Congruences relating to the sums of products of the first $n$ numbers and to other
sums of products\jour Quart. J. Math.\vol 31\yr 1900\pages 1--35\endref

\ref\key G2\by J.W.L. Glaisher\paper On the residues of the sums of products of the first $p-1$ numbers, and their powers, to modulus $p^2$ or $p^3$
\jour Quart. J. Math.\vol 31\yr 1900\pages 321--353\endref

\ref\key HT\by C. Helou and G. Terjanian\paper {\rm On Wolstenholme's theorem and its converse}
\jour J. Number Theory 128(2008) 475--499\endref

\ref\key IR\by K. Ireland and M. Rosen \book {\rm A Classical
Introduction to Modern Number Theory, second ed., Graduate Texts in
Math., Vol. 84,} \publ Springer, New York, 1990\endref


\ref\key L\by L. Long \paper Hypergeometric evaluation identities
and supercongruences \jour Pacific J. Math.\vol 249\yr 2011\pages
405--418\endref

\ref\key M\by R. Me\v strovi\'c\paper Proof of a congruence for harmonic numbers conjectured by Z.-W. Sun
\jour Int. J. Number Theory\vol 8\yr 2012\pages 1081--1085\endref

\ref\key S\by Z.-H. Sun\paper Congruences concerning Bernoulli numbers and Bernoulli polynomials
\jour Discrete Appl. Math.\vol 105\yr 2000\pages 193--223\endref

\ref\key Su1\by Z.-W. Sun\paper On sums of Apery polynomials and related congruences
 \jour J. Number Theory\vol 132\yr 2012\pages 2673--2699\endref

\ref\key Su2\by Z.-W. Sun\paper Arithmetic theory of harmonic
numbers \jour Proc. Amer. Math. Soc.\vol 140\yr 2012\pages 415--428 \endref

\ref\key Su3\by Z.-W. Sun\paper Conjectures and results on $x^2$ mod $p^2$ with $4p=x^2+dy^2$
 \jour in: Number Theory and Related Area (eds., Y. Ouyang, C. Xing, F. Xu and P. Zhang), Adv. Lect. Math. 27,
 Higher Education Press and Internat. Press, Beijing-Boston, 2013, pp. 149--197
\endref

\ref\key ST\by Z.-W. Sun and R. Tauraso\paper New congruences for central binomial coefficients
\jour Adv. in Appl. Math.\vol 45\yr 2010\pages 125--148\endref

\ref\key T\by R. Tauraso\paper More congruences for central binomial coefficients
\jour J. Number Theory \vol 130\yr 2010\pages 2639--2649\endref

\ref\key vLW\by J.H. van Lint and R. M. Wilson\book A Course in Combinatorics\ed 2nd\publ Cambridge Univ. Press, Cambridge, 2001\endref

\ref\key W\by J. Wolstenholme\paper On certain properties of prime numbers\jour Quart. J. Appl. Math.
\vol 5\yr 1862\pages 35--39\endref


\endRefs
\enddocument